\documentclass[11pt, reqno]{amsart}
\usepackage{amsmath,amsopn,amssymb,amsthm}
\usepackage[breaklinks=true,colorlinks=true,linkcolor=blue,citecolor=blue,urlcolor=blue]{hyperref}
\usepackage[cp1251]{inputenc}
\usepackage[english]{babel}

\DeclareMathOperator{\Ad}{Ad}
\DeclareMathOperator{\ad}{ad}

\begin{document}

\title{Addendum to ``Bounded Isometries and Homogeneous Quotients'',
JGA 27 (2017), 56--64}

\author{Joseph A. Wolf}

\maketitle

A few years ago Profs. Yurii Nikonorov and Ming Xu found a gap \cite{NX2019} 
in my paper
\cite{Wolf16} on homogeneity and bounded isometries.  In lines 2 and 3 of
the third paragraph of the proof of \cite[Theorem 2.5]{Wolf16} 
it was only shown that $Ad(sk)$ is bounded, and from that it should have been
shown that $Ad(s)$ is bounded.  Now Prof. Nikonorov 
modified my original proof, completing the argument, 
and agreed to its publication as a correction.  See below.
The statement of the theorem is unchanged; 
the modified argument uses a result \cite{Mos78} of Moskowitz.

\medskip
\noindent {\bf Theorem} (\cite[Theorem 2.5]{Wolf16}).
{\em Let $(M,d)$ be a metric space on which an exponential solvable Lie group 
$S$ acts effectively and transitively by isometries. Let $G = I(M,d)$. Then 
$G$ is a Lie group, any isotropy subgroup $K$ is compact, and $G = SK$. 
If $g\in G$ is a bounded isometry then $g$ is a central element in $S$.}

{\bf Proof} (Yurii Nikonorov \cite{N2021}).
$M$ carries a differentiable manifold structure for which $s \mapsto s(x_0)$ 
is a diffeomorphism $S\cong M$. In the compact--open topology $G = I(M, d)$ 
is locally compact and its its action on M is proper \cite{DW28}. In 
particular, if $x_0 \in M$ then the isotropy subgroup 
$K = \{k\in G\, \mid \, k(x_0) = x_0\}$ is compact. Further 
\cite[Corollary in \S 6.3]{MZ55} $G$ is a Lie group. Now $S$ and $K$ are 
closed subgroups, $G = SK$, and $M = G/K$. The action of $G$ on
$M = SK/K = S$ is $(sk) : s' \mapsto s\cdot k s' k^{-1}$.  We now assume 
that $x_0=1\in S$.

Express $g = sk \in G$ with $s\in S$ and $k\in K$. Suppose that $g$ is a 
bounded isometry of $(M, d)$. So there is a compact set $C\subset S$ such 
that $\tau(g)(s') s'^{-1}= s[(ks'k^{-1})s'^{-1}]\in C$ for every $s'\in S$.
Then $(ks'k^{-1})s'^{-1} \in s^{-1}C$ so the automorphism 
$a: S \rightarrow S$, $a(s')=ks'k^{-1}$
is an automorphism of bounded displacement. By Corollary 1.3 in~\cite{Mos78}, 
we see that $a$ is trivial. Since the isotropy representation is faithful,
it follows that $k=1$ and $g=s\in S$.

Express $g = s = \exp(\xi)$ where 
$\xi \in \mathfrak{s}:=\operatorname{Lie}(S)$. Decompose $\mathfrak{s}$ as a 
vector space direct sum $\mathfrak{n}+\mathfrak{a}$,
where $\mathfrak{n}$ is the nilradical of $\mathfrak{s}$.
In a basis respecting that direct sum, $\ad (\xi)|_{\mathfrak{s}}$ has 
matrix of the form
$\left( \begin{smallmatrix} 0 & h \\ 0 & 0 \\ \end{smallmatrix} \right)$. 
Let $N$ denote the unipotent radical of $S$, so 
$\mathfrak{n}:=\operatorname{Lie}(N)$.  Then
the $(1, 1)$ and $(2, 1)$ blocks vanish because $g$
centralizes $N$ by \cite[Th\'{e}or\`{e}me 1]{Tits64},
and the $(2, 2)$ block vanishes because 
$[\mathfrak{s},\mathfrak{s}]\subset \mathfrak{n}$.  Note
$\left( \begin{smallmatrix} 0 & h \\ 0 & 0 \\ \end{smallmatrix} \right)^2 = 0$, 
so $\Ad(g)|_{\mathfrak{s}}- I = \ad (\xi)|_{\mathfrak{s}}$. It cannot have 
relatively compact image unless $h = 0$. Thus $\Ad(g)|_{\mathfrak{s}}= I$, and
$g\in  (Z_G(S) \cap S)=Z_S$ as asserted.
\hfill $\square$

\bibliographystyle{amsunsrt}

\begin{thebibliography}{[99]}


\bibitem{DW28}
D. van Dantzig, B. L. van der Waerden, \"{U}ber metrisch homogene R\"{a}ume, Abh. Math. Seminar
Hamburg 6 (1928), 367--376.

\bibitem{MZ55}
D. Montgomery, L. Zippin, ``Topological Transformation Groups'', Interscience, 1955.

\bibitem{Mos78}
M.~Moskowitz,
Some remarks on automorphisms of bounded displacement and bounded cocycles,
Monatsh. Math. 85  (1978), 323--336.

\bibitem{N2021}
Y. Nikonorov,
private communication, 2021.

\bibitem{NX2019}
Y. Nikonorov \& Ming Xu,
private communications, 2019.

\bibitem{Tits64}
J. Tits,  Automorphismes \`a deplacement borne des groupes de Lie, Topology 3 (1964), 97--107.

\bibitem{Wolf16}
J. A. Wolf, Bounded isometries and homogeneous quotients, J. Geometric Anal. 27 (2017) 56--64.

\end{thebibliography}

\vspace{5mm}

\end{document}